\documentclass[english,11pt,reqno]{amsart}
\usepackage{amssymb,amsmath,amstext,amsthm,amsfonts}
\usepackage{graphicx}
\usepackage[ansinew]{inputenc}

\newcommand{\length}{\operatorname{length}}
\newcommand{\leb}{\operatorname{Leb}}

\newcommand{\dist}{\operatorname{dist}}
\newcommand{\diam}{\operatorname{diam}}

\begin{document}

\newcommand{\mcup}{\mbox{$\bigcup$}}
\newcommand{\mcap}{\mbox{$\bigcap$}}

\def \RR {{\mathbb R}}
\def \ZZ {{\mathbb Z}}
\def \NN {{\mathbb N}}
\def \PP {{\mathbb P}}
\def \TT {{\mathbb T}}
\def \II {{\mathbb I}}
\def \JJ {{\mathbb J}}

\def \vare {\varepsilon }

 \def \cf {\mathcal{F}}
 \def \cm {\mathcal{M}}
 \def \cn {\mathcal{N}}
 \def \cq {\mathcal{Q}}
 \def \cp {\mathcal{P}}
 \def \cb {\mathcal{B}}
 \def \cc {\mathcal{C}}
 \def \cs {\mathcal{S}}
 \def \bc {\mathcal{B}}
 \def \hc {\mathcal{C}}

\newcommand{\dem}{\begin{proof}}
\newcommand{\cqd}{\end{proof}}

\newcommand{\qand}{\quad\text{and}\quad}

\newtheorem{theorem}{Theorem}
\newtheorem{corollary}{Corollary}

\newtheorem*{Maintheorem}{Main Theorem}

\newtheorem{maintheorem}{Theorem}
\renewcommand{\themaintheorem}{\Alph{maintheorem}}
\newcommand{\cmt}{\begin{maintheorem}}
\newcommand{\fmt}{\end{maintheorem}}

\newtheorem{maincorollary}[maintheorem]{Corollary}
\renewcommand{\themaintheorem}{\Alph{maintheorem}}
\newcommand{\cmc}{\begin{maincorollary}}
\newcommand{\fmc}{\end{maincorollary}}

\newtheorem{T}{Theorem}[section]
\newcommand{\ct}{\begin{T}}
\newcommand{\ft}{\end{T}}

\newtheorem{Corollary}[T]{Corollary}
\newcommand{\cco}{\begin{Corollary}}
\newcommand{\fco}{\end{Corollary}}

\newtheorem{Proposition}[T]{Proposition}
\newcommand{\cpr}{\begin{Proposition}}
\newcommand{\fpr}{\end{Proposition}}

\newtheorem{Lemma}[T]{Lemma}
\newcommand{\cle}{\begin{Lemma}}
\newcommand{\fle}{\end{Lemma}}

\newtheorem{Sublemma}[T]{Sublemma}
\newcommand{\csle}{\begin{Lemma}}
\newcommand{\fsle}{\end{Lemma}}

\newtheorem{Remark}[T]{Remark}
\newcommand{\cre}{\begin{Remark}}
\newcommand{\fre}{\end{Remark}}

\newtheorem{Definition}[T]{Definition}
\newcommand{\cd}{\begin{Definition}}
\newcommand{\fd}{\end{Definition}}

\title[Partially hyperbolic horseshoes]{On the nonexistence of fat\\ partially hyperbolic horseshoes}

\author{Jos\'e F. Alves}
\address{Departamento de Matem\'atica Pura, Faculdade de Ci\^encias do Porto\\
Rua do Campo Alegre 687, 4169-007 Porto, Portugal}
\email{jfalves@fc.up.pt} \urladdr{http://www.fc.up.pt/cmup/jfalves}

\author{Vilton  Pinheiro}
\address{Departamento de Matem\'atica, Universidade Federal da Bahia\\
Av. Ademar de Barros s/n, 40170-110 Salvador, Brazil.}
\email{viltonj@ufba.br}

\date{\today}

\thanks{Work carried out at the Federal University of
Bahia. JFA is partially supported by CMUP and by a grant of FCT. VP
is partially supported by PADCT/CNPq}

\begin{abstract}
We show that there are no partially hyperbolic horseshoes with
positive Lebesgue measure for diffeomorphisms whose class of
differentiability is higher than 1. This generalizes a classical
result by Bowen for uniformly hyperbolic horseshoes.
\end{abstract}

\maketitle

\setcounter{tocdepth}{1}

\tableofcontents


\section{Introduction}

Let $M$ be a Riemannian manifold. We use $\leb$ to denote a
normalized volume form defined on the Borel sets of $M$ that we call
{\em Lebesgue measure.} Given a submanifold $\gamma\subset M$  we
use $\leb_\gamma$ to denote the measure on $\gamma$ induced by the
restriction of the Riemannian structure to $\gamma$.

Let $f:M\to M$ be a $C^{1}$ diffeomorphism,  and let $\Lambda\subset
M$ be a compact {\em invariant set}, i.e.
$f(\Lambda)\subset\Lambda$. We say that $\Lambda$ is a {\em
hyperbolic set} if there is a $Df$-invariant splitting $T_\Lambda
M=E^{s}\oplus E^{u}$ of the tangent bundle restricted to~$\Lambda$
and a constant $\lambda<1$ such that (for some choice of a
Riemannian metric on $M$) for every $x\in \Lambda$
      $$\|Df \mid E^{s}_x\|<\lambda\qand \|Df^{-1} \mid E^{u}_x\|<\lambda.$$

We say that an embedded disk $\gamma\subset M$ is an {\em unstable
manifold}, or an {\em unstable disk}, if
$\dist(f^{-n}(x),f^{-n}(y))\to0$ exponentially fast as $n\to\infty$,
for every $x,y\in\gamma$. Similarly, $\gamma$ is called a {\em
stable manifold}, or a {\em stable disk}, if
$\dist(f^{n}(x),f^{n}(y))\to0$ exponentially fast as $n\to\infty$,
for every $x,y\in\gamma$. It is well-known that every point in a
hyperbolic set possesses a local stable manifold $W_{loc}^s(x)$ and
a local unstable manifold $W_{loc}^u(x)$ which are disks tangent to
$E^{s}_x$ and $E^{u}_x$ at $x$ respectively.

A hyperbolic set $\Lambda$ is said to be a {\em horseshoe} if local
stable and local unstable manifolds through points in $\Lambda$
intersect $\Lambda$ in a Cantor set. Horseshoes were introduced by
Smale and appear naturally when one unfolds a homoclinic tangency
associated to some hyperbolic periodic point of saddle type. It
follows from \cite[Theorem 4.11]{B}  that a $C^{1+\alpha}$
diffeomorphism  cannot have a fat hyperbolic horseshoe, i.e. a
hyperbolic horseshoe $\Lambda$ with $\leb(\Lambda)>0$; actually the
result in~\cite{B} is proved for basic sets. Nevertheless, here we
obtain a generalization of that result to a much more general
situation. Let us remark that fat hyperbolic horseshoes exist for
$C^1$ diffeomorphisms, as shown in \cite{Bo}.

%

%

We say that a compact invariant set $\Lambda$ has a {\em dominated
splitting\/} if there exists a continuous $Df$-invariant splitting
$T_\Lambda M=E^{cs}\oplus E^{cu}$ of the tangent bundle restricted
to $\Lambda$, and a constant $0<\lambda<1$ such that (for some
choice of a Riemannian metric on $M$) for every $x\in\Lambda$
$$
 \|Df \mid E^{cs}_x\|
\cdot \|Df^{-1} \mid E^{cu}_{f(x)}\| \le\lambda.
$$
We call $E^{cs}$ the {\em centre-stable bundle} and $E^{cu}$ the
{\em centre-unstable bundle}. We say that $f$ is {\em non-uniformly
expanding along the centre-unstable direction} for
 $x\in \Lambda$ if
\begin{equation}
\label{e.lyapunov}\tag{NUE} \liminf_{n\to+\infty} \frac{1}{n}
    \sum_{j=1}^{n} \log \|Df^{-1} \mid E^{cu}_{f^j(x)}\|<0.
\end{equation}
Condition \ref{e.lyapunov} means that the derivative has {\em
expanding behavior in average} over the orbit of $x$. This implies
that $x$ has $\dim(E^{cu})$ positive
Lypaunov exponents in the $E^{cu}_x$ direction. 
As shown in \cite[Theorem C]{AAS}, if condition~\ref{e.lyapunov}
holds for every point in a compact invariant set   $\Lambda$, then
$E^{cu}$ is necessarily {\em uniformly expanding in $\Lambda$}, i.e.
there is  $0<\lambda<1$ such that
 $$\|Df^{-1}\mid E^{cu}_{f(x)}\| \le\lambda,\quad\text{for every
 $x\in\Lambda$.}$$
 A class of diffeomorphisms with a dominated splitting
$TM=E^{cs}\oplus E^{cu}$ for which \ref{e.lyapunov} holds Lebesgue
almost everywhere in $M$ and $E^{cu}$ is not uniformly expanding can
be found in \cite[Appendix A]{ABV}.

 \cmt\label{t:disco1}
 Let \( f: M\to M \) be a  \( C^{1+\alpha} \)
diffeomorphism and let $\Lambda\subset M$ have a dominated
splitting. If there is $H\subset\Lambda$ with $\leb(H)>0$ such that
\ref{e.lyapunov} holds for every $x\in H$, then $\Lambda$ contains
some local unstable disk.
 \fmt

We say that a compact invariant set $\Lambda$ is {\em partially
hyperbolic} if it has a dominated splitting $T_\Lambda
M=E^{cs}\oplus E^{cu}$ for which $E^{cs}$ is {\em uniformly
contracting} or $E^{cu}$ is {\em uniformly expanding}, meaning that
there is $0<\lambda<1$ such that
$0<\lambda<1$ such that $ \|Df\mid E^{cs}_{x}\| \le\lambda$ for
every $x\in \Lambda$, or $ \|Df^{-1}\mid E^{cu}_{f(x)}\| \le\lambda$
for every $x\in \Lambda$.

 The next result is a direct consequence of Theorem~\ref{t:disco1},
whenever $E^{cu}$ is uniformly expanding. If, on the other hand,
$E^{cs}$ is uniformly contracting, then we just have to apply
Theorem~\ref{t:disco1} to $f^{-1}$.

 \cmc\label{c:disco1}
 Let \( f: M\to M \) be a  \( C^{1+\alpha} \)
diffeomorphism and let $\Lambda\subset M$ be a partially hyperbolic
set with $\leb(\Lambda)>0$.
\begin{enumerate}
  \item If $E^{cs}$ is uniformly contracting, then $\Lambda$ contains a local stable disk.
  \item If $E^{cu}$ is uniformly expanding, then $\Lambda$ contains a local unstable disk.
\end{enumerate}

 \fmc

In particular, \( C^{1+\alpha} \) diffeomorphisms cannot have
partially hyperbolic horseshoes with positive Lebesgue measure. The
same holds for partially hyperbolic sets intersecting a local stable
or a local unstable disk in a positive Lebesgue measure subset, as
Corollary~\ref{c:disco2} below shows.

 \cmt\label{t:disco2}
 Let \( f: M\to M \) be a  \( C^{1+\alpha} \)
diffeomorphism and let $\Lambda\subset M$ have a dominated
splitting. Assume that there is a local unstable disk $\gamma$ with
$\leb_\gamma(\gamma\cap\Lambda)>0$ such that \ref{e.lyapunov} holds
for every $x\in\gamma\cap\Lambda$. Then $\Lambda$ contains some
local unstable disk.
 \fmt

The next result is a direct consequence of Theorem~\ref{t:disco2} in
the case that $E^{cu}$ is uniformly expanding, and a consequence of
the theorem applied to $f^{-1}$ in the case that $E^{cs}$ is
uniformly contracting.

 \cmc\label{c:disco2}
 Let \( f: M\to M \) be a  \( C^{1+\alpha} \)
diffeomorphism and let $\Lambda\subset M$ be a partially hyperbolic
set.
\begin{enumerate}
\item If $E^{cs}$ is uniformly contracting and there is a local stable disk
$\gamma$ such that $\leb_\gamma(\gamma\cap\Lambda)>0$, then
$\Lambda$ contains a local stable disk.
\item If $E^{cu}$ is
uniformly expanding and there is a local unstable disk $\gamma$ such
that $\leb_\gamma(\gamma\cap\Lambda)>0$, then $\Lambda$ contains a
local unstable disk.
\end{enumerate}

 \fmc

 Theorems~\ref{t:disco1} and \ref{t:disco2} are in fact corollaries
of a slightly more general result that we present at the beginning
of Section~\ref{s.local}.

\subsection*{Acknowledgement} We are grateful to M. Viana for
valuable references on this topic.

\section{H\"older control of tangent direction} This section is
a survey of results in~\cite[Section 2]{ABV} concerning the H\"older
control of the tangent direction of submanifolds. As observed in
\cite[Remark 2.3]{ABV} those results are valid for diffeomorphisms
of class $C^{1+\alpha}$. In this section we only use the existence
of a dominated splitting $T_\Lambda M=E^{cs} \oplus E^{cu}$. We fix
continuous extensions of the two bundles $E^{cs}$ and $E^{cu}$ to
some neighborhood $U$ of $\Lambda$, that we denote by
$\tilde{E}^{cs}$ and $\tilde{E}^{cu}$. We do not require these
extensions to be invariant under $Df$. Given $0<a<1$, we define the
{\em centre-unstable cone field
$C_a^{cu}=\left(C_a^{cu}(x)\right)_{x\in U}$ of width $a$\/} by
\begin{equation}
\label{e.cucone} C_a^{cu}(x)=\big\{v_1+v_2 \in
\tilde{E}_x^{cs}\oplus \tilde{E}_x^{cu} \text{\ such\ that\ }
\|v_1\| \le a \|v_2\|\big\}.
\end{equation}
We define the {\em centre-stable cone field
$C_a^{cs}=\left(C_a^{cs}(x)\right)_{x\in U}$ of width $a$\/} in a
similar way, just reversing the roles of the subbundles in
(\ref{e.cucone}). We fix $a>0$ and $U$ small enough so that, up to
slightly increasing $\lambda<1$, the domination condition remains
valid for any pair of vectors in the two cone fields:
$$
\|Df(x)v^{cs}\|\cdot\|Df^{-1}(f(x))v^{cu}\|
\le\lambda\|v^{cs}\|\,\|v^{cu}\|
$$
for every $v^{cs}\in C_a^{cs}(x)$, $v^{cu}\in C_a^{cu}(f(x))$, and
any point $x\in U\cap f^{-1}(U)$. Note that the centre-unstable cone
field is positively invariant: $$Df(x) C_a^{cu}(x)\subset
C_a^{cu}(f(x)),\quad\text{whenever $x,f(x)\in U$.}$$ Indeed, the
domination property together with the invariance of
$E^{cu}=(\tilde{E}^{cu} \mid \Lambda)$ imply that
$$
Df(x) C_a^{cu}(x) \subset C_{\lambda a}^{cu}(f(x))
                  \subset C_a^{cu}(f(x)),
$$
for every $x\in K$. This extends to any $x\in U\cap f^{-1}(U)$ just
by continuity.



We say that an embedded $C^1$ submanifold $N\subset U$ is {\em
tangent to the centre-unstable cone field\/} if the tangent subspace
to $N$ at each point $x\in N$ is contained in the corresponding cone
$C_a^{cu}(x)$. Then $f(N)$ is also tangent to the centre-unstable
cone field, if it is contained in $U$, by the domination property.

Our aim now is to express the notion of H\"older variation of the
tangent bundle in local coordinates. We choose $\delta_0>0$ small
enough so that the inverse of the exponential map $\exp_x$ is
defined on the $\delta_0$ neighbourhood of every point $x$ in $U$.
From now on we identify this neighbourhood of $x$ with the
corresponding neighbourhood $U_x$ of the origin in $T_x N$, through
the local chart defined by $\exp_x^{-1}$. 
Reducing $\delta_0$, if necessary, we may suppose that
$\tilde{E}^{cs}_x$ is contained in the centre-stable cone
$C^{cs}_a(y)$ of every $y\in U_x$. In particular, the intersection
of $C^{cu}_a(y)$ with $\tilde{E}^{cs}_x$ reduces to the zero vector.
Then, the tangent space to $N$ at $y$ is parallel to the graph of a
unique linear map $A_x(y):T_x N \to \tilde{E}_x^{cs}$. Given
constants $C>0$ and $0<\zeta\le 1$, we say that {\em the tangent
bundle to $N$ is $(C,\zeta)$-H\"older\/} if for every $y\in N \cap
U_x $ and $x\in V_0$
\begin{equation}
\label{e.holder} \|A_x(y)\|\le C d_x(y)^\zeta ,
\end{equation}
where $d_x(y)$ denotes the distance from $x$ to $y$ along $N\cap
U_x$, defined as the length of the shortest curve connecting $x$ to
$y$ inside $N\cap U_x$.

Recall that we have chosen the neighbourhood $U$ and the cone width
$a$ sufficiently small so that the domination property remains valid
for vectors in the cones $C_a^{cs}(z)$, $C_a^{cu}(z)$, and for any
point $z$ in $U$. Then, there exist $\lambda_1 \in (\lambda,1)$ and
$\zeta\in(0,1]$ such that
\begin{equation}
\label{e.dominacao} \|Df(z) v^{cs}\| \cdot \|Df^{-1}(f(z))
v^{cu}\|^{1+\zeta} \le \lambda_1 < 1
\end{equation}
for every norm $1$ vectors $v^{cs}\in C_a^{cs}(z)$ and $v^{cu}\in
C_a^{cu}(z)$, at any $z\in U$. Then, up to reducing $\delta_0>0$ and
slightly increasing $\lambda_1<1$, condition (\ref{e.dominacao})
remains true if we replace $z$ by any $y\in U_x$, $x\in U$ (taking
$\|\cdot\|$ to mean the Riemannian metric in the corresponding local
chart).

We fix $\zeta$ and $\lambda_1$ as above. Given a $C^1$ submanifold
$N\subset U$, we define
\begin{equation}
\label{e.kappa} \kappa(N)=\inf\{C>0:\text{the tangent bundle of $N$
is $(C,\zeta)$-H\"older}\}.
\end{equation}
The next result appears in \cite[Corollary 2.4]{ABV}.

%
%

\cpr \label{c.curvature} There exists $C_1>0$ such that, given any
$C^1$ submanifold $N\subset U$ tangent to the centre-unstable cone
field,
\begin{enumerate}
\item there exists $n_0\ge 1$ such that $\kappa(f^n(N)) \le
C_1$ for every $n\ge n_0$ such that $f^k(N) \subset U$ for all $0\le
k \le n$;
 \item if $\kappa(N) \le C_1$, then the same is true
for every iterate $f^n(N)$  such that $f^k(N)\subset U$ for all
$0\le k \le n$;
 \item in particular, if $N$ and $n$ are
as in (2), then the functions
$$
J_k: f^k(N)\ni x \mapsto \log |\det \big(Df \mid T_x f^k(N)\big)|,
\quad\text{$0\le k \le n$},
$$
are $(L,\zeta)$-H\"older continuous with $L>0$ depending only on
$C_1$ and~$f$.
\end{enumerate}
\fpr

\section{Hyperbolic times and bounded distortion}\label{sub.hyptimes}

The following notion will allow us to derive {\em uniform
behaviour\/} (expansion, distortion) from the non-uniform expansion.

\cd \label{d.hyperbolic1} Given $\sigma<1$, we say that $n$ is a
{\em $\sigma$-hyperbolic time\/} for  $x\in \Lambda$ if
$$
\prod_{j=n-k+1}^{n}\|Df^{-1} \mid E^{cu}_{f^{j}(x)}\| \le \sigma^k,
\quad\text{for all $1\le k \le n$.}
$$
\fd In particular, if $n$ is a $\sigma$-hyperbolic time for $x$,
then $Df^{-k} \mid E^{cu}_{f^{n}(x)}$ is a contraction for every
$1\le k \le n$:
$$
\|Df^{-k} \mid E^{cu}_{f^{n}(x)}\| \le \prod_{j=n-k+1}^{n}\|Df^{-1}
\mid E^{cu}_{f^{j}(x)}\| \le \sigma^{k}.
$$
Moreover, if $a>0$ is taken sufficiently small in the definition of
our cone fields, and we choose $\delta_1>0$ also small so that the
$\delta_1$-neighborhood of $\Lambda$ should be contained in $U$,
then by continuity
\begin{equation}
\label{e.delta1} \|Df^{-1}(f(y)) v \| \le \frac{1}{\sqrt\sigma}
\|Df^{-1}|E^{cu}_{f(x)}\|\,\|v\|,
\end{equation}
whenever $x\in \Lambda$, $\dist(x,y)\le \delta_1$, and $v\in
C^{cu}_a(y)$.


Given any disk $\Delta\subset M$, we use $\dist_\Delta(x,y)$ to
denote the distance between $x,y\in \Delta$, measured along
$\Delta$. The distance from a point $x\in \Delta$ to the boundary of
$\Delta$ is $\dist_\Delta(x,\partial \Delta)= \inf_{y\in\partial
\Delta}\dist_\Delta(x,y)$.

\smallskip

\cle \label{l.contraction} Take any $C^1$ disk $\Delta\subset U$ of
radius $\delta$, with $0<\delta<\delta_1$, tangent to the
centre-unstable cone field.
There is $n_0\ge1$ such that if $x\in \Delta$ with
$\dist_\Delta(x,\partial \Delta)\ge \delta/2$ and $n \ge n_0$  is a
$\sigma$-hyperbolic time for $x$, then there is a neighborhood $V_n$
of $x$ in $\Delta$ such that:
\begin{enumerate}
    \item $f^{n}$ maps $V_n$ diffeomorphically onto a
disk  of radius $\delta_1$ around  $f^{n}(x)$ tangent to the
centre-unstable cone field;
    \item for every $1\le k
\le n$ and $y, z\in V_n$, $$
\dist_{f^{n-k}(V_n)}(f^{n-k}(y),f^{n-k}(z)) \le
\sigma^{k/2}\dist_{f^n(V_n)}(f^{n}(y),f^{n}(z)).$$
\end{enumerate}
\fle

\begin{proof}
First we show that $f^n(\Delta)$ contains some disk of radius
$\delta_1$ around $f^n(x)$, as long as
 \begin{equation}\label{e.ene}
 n>2\frac{\log(\delta/(2\delta_1))}{\log(\sigma)}.
 \end{equation}
 Assume that there is $y\in \partial \Delta$ with
$\dist_{f^n(\Delta)}(f^n(x),f^n(y))<\delta_1$. Let $\eta_0$ be a
curve of minimal length in $f^n(\Delta)$ connecting $f^n(x)$ to
$f^n(y)$. For $0\le k \le n$ we write $\eta_k=f^{n-k}(\eta_0)$. We
prove  by induction that $\length(\eta_k) < \sigma^{k/2}\delta_1, $
for $0 \le k \le n$. Let $1 \le k \le n$ and assume that
$$
\length(\eta_j) < \sigma^{j/2}\delta_1, \quad\text{for\ }0 \le j \le
k-1.
$$
Denote by $\dot\eta_0(w)$ the tangent vector to the curve $\eta_0$
at the point $w$. Then, by the choice of $\delta_1$ in
(\ref{e.delta1}) and the definition of $\sigma$-hyperbolic time,
$$
\|Df^{-k}(w) \dot\eta_0(w)\| \le {\sigma^{-k/2} \,
\|\dot\eta_0(w)\|}
         {\prod_{j=n-k+1}^{n}\|Df^{-1} |E^{cu}_{f^j(x)}\|}
\le \sigma^{k/2}\|\dot\eta_0(w)\|.
$$
Hence,
$$
\length(\eta_k) \le \sigma^{k/2}\length(\eta_0) < \sigma^{k/2}
\delta_1\,.
$$
This completes our induction. In particular, we have $
\length(\eta_n)  < \sigma^{n/2} \delta_1\,. $ Note that $\eta_n$ is
a curve in $\Delta$ connecting $x$ to $y\in
\partial \Delta$, and so $\length(\eta_n)\ge
\delta/2$. Thus we must have
 \begin{equation*}\label{e.ene2}
 n<2\frac{\log(\delta/(2\delta_1))}{\log(\sigma)}.
 \end{equation*}
 Hence  $f^n(\Delta)$ contains some disk of radius
$\delta_1$ around $f^n(x)$ for $n$ as in~\eqref{e.ene}.

Let now $\Delta_1$ be the disk of radius $\delta_1$ around $f^n(x)$
in $f^n(\Delta)$ and let $V_n=f^{-n}(\Delta_1)$,  for $n$ as in
\eqref{e.ene}. Take any $y,z\in V_n$ and let $\eta_0$ be a curve of
minimal length in $\Delta_1$ connecting $f^n(y)$ to $f^n(z)$.
Defining $\eta_k=f^{n-k}(\eta_0)$, for $1\le k \le n$, and arguing
as before we inductively prove that for $1\le k\le n$
$$
\length(\eta_k) \le \sigma^{k/2}\length(\eta_0) =
\sigma^{k/2}\dist_{f^{n}(V_n)}(f^{n}(y),f^{n}(z)) ,
$$
which implies that for $1\le k \le n$
 $$
\dist_{f^{n-k}(V_n)}(f^{n-k}(y),f^{n-k}(z)) \le
\sigma^{k/2}\dist_{f^n(V_n)}(f^{n}(y),f^{n}(z)).$$
This completes the proof of the lemma.
\end{proof}

We shall sometimes refer to the sets \( V_n \) as \emph{hyperbolic
pre-balls} and to their images \( f^{n}(V_n) \) as \emph{hyperbolic
balls}. Notice that the latter are indeed balls of radius \(
\delta_1 \).

\cco[\bf Bounded Distortion] \label{p.distortion} There exists
$C_2>1$ such that given $\Delta$ as in Lemma~\ref{l.contraction}
with $\kappa(\Delta) \le C_1$, and given any hyperbolic pre-ball
$V_n\subset \Delta$ with $n\ge n_0$, then for all $y,z\in V_n$
%
$$
\frac{1}{C_2} \le \frac{|\det Df^{n} \mid T_y \Delta|}
                     {|\det Df^{n} \mid T_z \Delta|}
            \le C_2.
$$
 \fco

\begin{proof}
For $0\le i <n$ and $y\in \Delta$, we denote $J_i(y)= |\det Df \mid
T_{f^i(y)} f^i(\Delta)|$. Then,
$$
\log \frac{|\det Df^{n} \mid T_y \Delta|}{|\det Df^{n} \mid T_z
\Delta|} = \sum_{i=0}^{n-1} \big(\log J_i(y) - \log J_i(z)\big).
$$
By Proposition~\ref{c.curvature}, $\log J_i$ is $(L,\zeta)$-H\"older
continuous, for some uniform constant $L>0$. Moreover, by
Lemma~\ref{l.contraction}, the sum of all $\dist_\Delta(f^j(y),
f^j(z))^\zeta$ over $0\le j \le n$ is bounded by
$\delta_1/(1-\sigma^{\zeta/2})$. Now it suffices to take
$C_2=\exp(L\delta_1/(1-\sigma^{\zeta/2}))$.
\end{proof}


\section{A local unstable disk inside $\Lambda$}\label{s.local}

Now we are able to prove Theorems~\ref{t:disco1} and~\ref{t:disco2}.
These will be obtained as corollaries of the next result, as we
shall see below.

 \ct\label{t:disco}
 Let \( f: M\to M \) be a  \( C^{1+\alpha} \)
diffeomorphism and let $\Lambda\subset M$ have a dominated
splitting. Assume that there is a disk $\Delta$ tangent to the
centre-unstable cone field with $\leb_\Delta(\Delta\cap\Lambda)>0$
such that \ref{e.lyapunov} holds for every $x\in \Delta\cap\Lambda$.
Then $\Lambda$ contains some local unstable disk.
 \ft

Assume that there is $H\subset\Lambda$ with $\leb(H)>0$ such that
\ref{e.lyapunov} holds for every $x\in H$. Choosing a $\leb$ density
point of $H$, we laminate a neighborhood of that point into disks
tangent to the centre-unstable cone field contained in $U$. Since
the relative Lebesgue measure of the intersections of these disks
with $H$ cannot be all equal to zero, we obtain some disk $\Delta$
as in Theorem~\ref{t:disco} under the assumption of
Theorem~\ref{t:disco1}. For Theorem~\ref{t:disco2}, observe that
local stable manifolds are tangent to the centre-unstable spaces and
these vary continuously with the points in $\Lambda$, thus being
locally tangent to the centre-unstable cone field.

Let us now prove Theorem~\ref{t:disco}. Let $\Delta$ be a disk
tangent to the centre-unstable cone field intersecting $\Lambda$ in
a positive $\leb_\Delta$ subset such that \ref{e.lyapunov} holds for
every $x\in \Delta\cap\Lambda$. Let $H=\Delta\cap\Lambda$. Taking a
subset of $H$, if necessary,  still with positive $\leb_\Delta$
measure, we may assume that there is $c>0$ such that for every $x\in
H$
\begin{equation}
\label{e.lyapunovc} \liminf_{n\to+\infty} \frac{1}{n}
    \sum_{j=1}^{n} \log \|Df^{-1} \mid E^{cu}_{f^j(x)}\|\le-c.
\end{equation}
Since condition \eqref{e.lyapunovc} remains valid under iteration,
by Proposition~\ref{c.curvature} we may assume that
$\kappa(\Delta)<C_1$. It is no restriction to assume that $H$
intersects the sub-disk of $\Delta$ of radius $\delta/2$, for some
$0<\delta<\delta_1$, in a positive $\leb_\Delta$ subset, and we do
so.

The following lemma is due to Pliss \cite{Pl72}, and a proof of it
in the present form can be found in \cite[Lemma 3.1]{ABV}.

\cle[\bf Pliss] \label{l.pliss} Given $A \ge c_2>c_1>0$, let
$\theta=(c_2-c_1)/(A-c_1)$. Given any real numbers $a_1, \ldots,
a_N$ such that $a_j \le A$ and
$$
\sum_{j=1}^{N} a_j \ge c_2 N,\quad \text{for every } 1\le j \le N,
$$
then there are $l > \theta N$ and $1 < n_1 < \cdots < n_l \le N$ so
that
$$
\sum_{j=n+1}^{n_i} a_j\ge c_1 (n_i-n), \quad\text{for every } 0 \le
n < n_i \text{ and  }i=1, \dots, l.
$$
\fle

\cco\label{l:hyperbolic2}
    There is $\sigma>0$ 
     such that every \( x
    \in  H \) has infinitely many $\sigma$-hyperbolic times.
     \fco
\dem Given $x\in H$, by \eqref{e.lyapunovc} we have infinitely many
values  $N$ for which
$$
\sum_{j=1}^{N} \log \|Df^{-1}|E^{cu}_{f^j(x)}\| \le -\frac{c}2N.
$$
 Then it suffices to take
$c_1=c/2$, $c_2=c$, $A=\sup\big|\,\log \|Df^{-1}|E^{cu}\|\,\big|$,
and $a_j=-\log \|Df^{-1}|E^{cu}_{f^j(x)}\|$ in the previous lemma.
\cqd

Note that under assumption \eqref{e.lyapunovc} we are unable to
prove the existence of a positive frequency of hyperbolic times at
infinity. This would be possible if we had $\limsup$ instead of
$\liminf$ in \eqref{e.lyapunovc}, as shown in \cite[Corollary
3.2]{ABV}. The existence of infinitely many hyperbolic times is
enough for what comes next.

\cle \label{l.discao} There are an infinite sequence of integers
$1\le k_1<k_2<\cdots$ and, for each $n\in\NN$,  a disk $\Delta_n$ of
radius $\delta_1/4$ tangent to the centre-unstable cone field   such
that the relative Lebesgue measure of the set $f^{k_n}(H)$ in
$\Delta_n$ converges to 1 as $n\to\infty$.
\fle

\dem
%

Let $\epsilon>0$ be some small number. Let  $K$ be a compact subset
of $H$ and $A$ be an open neighborhood of $H$ in $\Delta$ such that
$$\leb_\Delta(A\setminus K)<\epsilon{\leb_\Delta(K)}.$$
It follows from Corollary~\ref{l:hyperbolic2} and
Lemma~\ref{l.contraction}  that we can choose for each $x\in K$ a
$\sigma$-hyperbolic time $n(x)$ and a hyperbolic pre-ball $V_x$ such
that $V_x\subset A$. Here $V_x$ is the neighborhood of $x$
associated to the hyperbolic time $n(x)$ constructed in
Lemma~\ref{l.contraction}, which is mapped diffeomorphically by
$f^{n(x)}$ onto a ball $B_{\delta_1}(f^{n(x)}(x))$ of radius
$\delta_1$ around $f^{n(x)}(x)$ tangent to the centre-unstable cone
field. Let $W_{x}\subset V_x$ be the pre-image of the ball
$B_{\delta_1/4}(f^{n(x)}(x))$ of radius $\delta_1/4$ under this
diffeomorphism.


By compactness we have $K\subset W_{x_1}\cup...\cup W_{x_s}$ for
some $x_1,...,x_m\in K$. Writing
\begin{equation}\label{eqn1}
    \{n_1,...,n_s\}=\{n(x_1),...,n(x_m)\},
\quad\text{with $n_1<n_2<...<n_s$},
\end{equation}
let $\mathcal{U}_1\subset\mathbb{N}$ be a maximal set of
$\{1,...,m\}$ such that if $u\in\mathcal{U}_1$ then $n(x_u)=n_1$ and
$W_{x_u}\cap W_{x_a}=\emptyset$ for all $a \in\mathcal{U}_1$ with
$a\ne u$. Inductively we define $\mathcal{U}_j$  for $1<j\le s$ as
follows. Suppose that $\mathcal{U}_{j-1}$ has already been defined.
Let $\mathcal{U}_j\subset\mathbb{N}$ be a maximal set of
$\{1,\dots,m\}$ such that if $u\in\mathcal{U}_j$ then $n(x_u)=n_j$
and $W_{x_u}\cap W_{x_a}=\emptyset$ for all $a \in\mathcal{U}_j$
with $a\neq u$, and also $W_{x_u}\cap W_{x_a}=\emptyset$ for all $
a\in \mathcal{U}_1\cup...\cup\mathcal{U}_{j-1}$.

Let $\mathcal{U}=\mathcal{U}_1\cup\cdots\cup\mathcal{U}_s$. By
maximality, each $W_{x_i}$, $1\le i\le m$, intersects some $W_{x_u}$
with $u\in\mathcal{U}$ and $n(x_i)\ge n(x_u)$. Thus, given any $1\le
i\le m$ and taking $u\in\mathcal{U}$ such that $W_{x_i}\cap
W_{x_u}\ne\emptyset$ and $n(x_i)\ge n(x_u)$, we get
$$
f^{n(x_u)}(W_{x_i})\cap B_{\delta_1/4}(f^{n(x_u)}(x_u))\ne\emptyset.
$$
Lemma~\ref{l.contraction}  assures that
$$
\diam(f^{n(x_u)}(W_{x_i}))\le\frac{\delta_1}{2}
\sigma^{(n(x_i)-n(x_u))/2}\le\frac{\delta_1}{2},
$$
 and so
$$f^{n(x_u)}(W_{x_i})\subset B_{\delta_1}(f^{n(x_u)}(x_u)).$$ This
implies that $W_{x_i}\subset V_{x_u}$. Hence
$\{V_{x_u}\}_{u\in\mathcal{U}}$ is a covering of $K$.

It follows from Corollary~\ref{p.distortion} that there is a uniform
constant $\gamma>0$ such that
 $$
 \frac{\leb_\Delta(W_{x_u})}{\leb_\Delta(V_{x_u})}\ge\gamma, \quad\text{for every }u\in\mathcal U.$$
Hence
 \begin{eqnarray*}
   \leb_\Delta\big(\cup_{u\in\mathcal{U}}W_{x_u}\big) &=& \sum_{u\in\mathcal{U}} \leb_\Delta(W_{x_u}) \\
    &\ge& \sum_{u\in\mathcal{U}}\gamma \leb_\Delta(V_{x_u}) \\
   &\ge& \gamma\leb_\Delta\big(\cup_{u\in\mathcal{U}}V_{x_u}\big)\\
   &\ge& \gamma \leb_\Delta(K).
 \end{eqnarray*}
Setting
$$
\rho=\min\left\{\frac{\leb_\Delta(W_{x_u}\setminus
K)}{\leb_\Delta(W_{x_u})}\colon u\in\mathcal{U}\right\},
$$
we have
\begin{eqnarray*}
  \varepsilon \leb_\Delta(K) &\ge& \leb_\Delta(A\setminus K) \\
   &\ge& \leb_\Delta\big(\cup_{u\in\mathcal{U}}W_{x_u}\setminus
K\big) \\
   &\ge& \sum_{u\in\mathcal{U}}\leb_\Delta(W_{x_u}\setminus K)\\
   &\ge& \rho \leb_\Delta\big(\cup_{u\in\mathcal{U}}W_{x_u}\big)\\
    &\ge& \rho \gamma
\leb_\Delta(K).
\end{eqnarray*}
 This implies that $\rho<\varepsilon/\gamma$. Since $\varepsilon>0$ can
 be taken arbitrarily small, 
%
we may always choose $W_{x_u}$ such that the relative Lebesgue
measure of $K$ in $W_{x_u}$ is arbitrarily close to 1. Then, by
bounded distortion, the relative Lebesgue measure of
$f^{n(x_u)}(H)\supset f^{n(x_u)}(K)$ in $f^{n(x_u)}(W_{x_u})$, which
is a disk of radius $\delta_1/4$ around $f^{n(x_u)}(x_u)$ tangent to
centre-unstable cone field, is also arbitrarily close to $1$.
Observe that since points in $H$ have infinitely many hyperbolic
times, we may take the integer $n(x_u)$ arbitrarily large, as long
as $n_1$ in \eqref{eqn1} is also taken large enough.
%
%
\cqd

\cpr\label{locdisk} There is a local unstable disk $\Delta_\infty$
of radius $\delta_1/4$ inside~$\Lambda$. \fpr

\dem Let $(\Delta_n)_n$ be the sequence of disks given by
Lemma~\ref{l.discao}, and consider $(x_n)_n$ the sequence of points
at which these disks are centered. Up to taking subsequences, we may
assume that the centers of the disks converge to some point $x$.
Using Ascoli-Arzela, the disks converge to some disk $\Delta_\infty$
centered at $x$. By construction, every point in $\Delta_\infty$ is
accumulated by some orbit of a point in $H\subset\Lambda$, and so
$\Delta_\infty\subset \Lambda$.

Note that each $\Delta_n$ is contained in the $k_n$-iterate of
$\Delta$, which is a disk tangent to the centre-unstable cone field.
The domination property implies that the angle between $\Delta_n$
and $E^{cu}$ goes to zero as $n\to\infty$, uniformly on $\Lambda$.
In particular, $\Delta_\infty$ is tangent to $E^{cu}$ at every point
in $\Delta_\infty\subset\Lambda$. By Lemma~\ref{l.contraction},
given any $k\ge 1$, then $f^{-k}$ is a $\sigma^{k/2}$-contraction on
$\Delta_n$ for every large $n$. Passing to the limit, we get that
$f^{-k}$ is a $\sigma^{k/2}$-contraction on $\Delta_\infty$ for
every $k\ge1$.

In particular, we have shown that the subspace $E_x^{cu}$ is
uniformly expanding for $Df$. The fact that the $Df$-invariant
splitting $T_\Lambda M=E^{cs}\oplus E^{cu}$ is a dominated splitting
implies  that any expansion $D f$ may exhibit along the
complementary direction $E^{cs}_x$ is weaker than the expansion in
the $E^{cu}_x$  direction. Then, by \cite{Pe76}, there exists a
unique unstable manifold $W^{u}_{loc}(x)$ tangent to $E^{cu}$ and
which is
contracted by the negative iterates of $f$. 
Since $\Delta_\infty$ is contracted by every $f^{-k}$, and  all its
negative iterates are tangent to centre-unstable cone field, then
$\Delta_\infty$ is contained in $W^{u}_{loc}(x)$.
 \cqd

\end{document}